\documentclass[a4paper,12pt]{article}
    \usepackage[top=2.5cm,bottom=2.5cm,left=2.5cm,right=2.5cm]{geometry}
    \usepackage{cite, amsmath, amssymb}
   \usepackage{tikz}
\begin{document}
\begin{center}
{\LARGE\bf Multiplicity of eigenvalues of cographs}
\end{center}
\begin{center}
{\large \bf Luiz Emilio Allem$^a,$ Fernando Tura$^{b,}$ \footnote{Corresponding author} }
\end{center}
\begin{center}
\it $^a$ Instituto  de Matem\'atica, UFRGS, Porto Alegre, RS, 91509-900, Brazil\\
\tt emilio.allem@ufrgs.br
\end{center}
\begin{center}
\it $^b$ Departamento de Matem\'atica, UFSM, Santa Maria, RS, 97105-900, Brazil\\
\tt ftura@smail.ufsm.br
\end{center}



\newcommand{\casei}{{\bf case~1}}
\newcommand{\subia}{{\bf subcase~1a}}
\newcommand{\subib}{{\bf subcase~1b}}
\newcommand{\subic}{{\bf subcase~1c}}
\newcommand{\caseii}{{\bf case~2}}
\newcommand{\subiia}{{\bf subcase~2a}}
\newcommand{\subiib}{{\bf subcase~2b}}
\newcommand{\caseiii}{{\bf case~3}}
\newcommand{\myvar}{x}
\newtheorem{Thr}{Theorem}
\newtheorem{Pro}{Proposition}
\newtheorem{Que}{Question}
\newtheorem{Con}{Conjecture}
\newtheorem{Cor}{Corollary}
\newtheorem{Lem}{Lemma}
\newtheorem{Fac}{Fact}
\newtheorem{Ex}{Example}
\newtheorem{Def}{Definition}
\newtheorem{Prop}{Proposition}
\def\floor#1{\left\lfloor{#1}\right\rfloor}

\newenvironment{my_enumerate}{
\begin{enumerate}
  \setlength{\baselineskip}{14pt}
  \setlength{\parskip}{0pt}
  \setlength{\parsep}{0pt}}{\end{enumerate}
}

\newenvironment{my_description}{
\begin{description}
  \setlength{\baselineskip}{14pt}
  \setlength{\parskip}{0pt}
  \setlength{\parsep}{0pt}}{\end{description}
}



\begin{abstract}
 Motivated by the linear time  algorithm that locates the eigenvalues of a cograph  $G$ \cite{JTT2016}, we  investigate the multiplicity of eigenvalue $\lambda$ for $ \lambda \neq 0,-1.$   For cographs with {\em balanced cotrees} we determine explicitly  the highest value for the multiplicity. The energy of a graph is defined as the sum of absolute values of the eigenvalues. A graph $G$ on $n$ vertices is said to be borderenergetic if its energy equals the energy of the complete graph $K_n.$ We present families of non-cospectral and borderenergetic cographs.
\end{abstract}

\baselineskip=0.30in

\section{Introduction}
\label{intro}
We recall that the \emph{spectrum} of a graph $G$ is the multiset of the eigenvalues of its adjacency matrix. The main goal of this paper is to discuss the multiplicity of eigenvalues of cographs.

Cographs is an important class of graphs for its many applications. They have several alternative characterizations, for example,  a cograph is graph which contains no  path  of length  four  as an induced subgraph \cite{Stewart} and because of this they are often simply  called  $P_4$ {\em free} graph in the literature. In particular it well known that any cograph has  a canonical tree representation, called  the {\em cotree} \cite{BSS2011}. The cotree will be relevant to this paper and will be described later.

Our original motivation for considering cographs is to study the distribution of eigenvalues of graphs. It is known, for example, that any interval of the real line contains some eigenvalues of graphs, since, more generally, any root of a real-rooted monic polynomial with integer coefficients occurs as an eigenvalue of some tree \cite{Sal2015}. On the other hand it was proved (see \cite{Moha}) that no cograph has eingenvalues in the interval $(-1,0)$, a surprising result.  

In this paper, we turn to study the multiplicities of eigenvalues of cographs. In \cite{JTT2013} it was proved that all eigenvalues of threshold graphs (a subclass of cographs), except $-1$ and $0$  are simple. This motivates us to investigate further the multiplicities of cograph eigenvalues. Since the multiplicities of the eigenvalues $-1$ and $0$ are known \cite{BSS2011} we deal with eigenvalues that are different from 0 and -1.

The multiplicities of graph eigenvalues are extensively studied by several authors. Bell {\em et al.} \cite{Bell}
determined upper bound for the multiplicities of graphs.  Later, Rowlinson in \cite{Row}  studied the multiplicities of eigenvalues in trees. Recently, Bu {\em et al.} \cite{Bu} studied the multiplicities in graphs attaching one pendent path,  generalizing some known results for trees and unicyclic graphs \cite{Row}.

Different from the star complement technique used in the works above,  our technique is based on an algorithm called {\em Diagonalization}, presented in \cite{JTT2016}.  The  Diagonalization  finds, in $O(n)$ time, the number of eigenvalues of a cograph, by operating directly on the cotree of the cograph. The algorithm and the technique will be explained in the next section.

We study cographs whose cotree is \emph{balanced} (see definition in Section \ref{multi}) and determine the multiplicity of some eigenvalues and an upper bound for the multiplicity of other eigenvalues. 

As an application of these results, we study the \emph{energy} of families of cographs. Recall that if $G$ is a graph having eigenvalues $\lambda_1, \ldots, \lambda_n,$ its  energy, denoted  $E(G)$  is defined  \cite{Gutman2015, Gutman2012} as $ \sum_{i=1}^n |  \lambda_i |.$ There are many results on energy and its applications in several areas, including in chemistry see  \cite{Gutman2012} for more details and the references therein.

It is  well known that the complete graph $K_n$ has  $E(K_n) = 2n-2$ and it is a natural and important research problem to determine graphs that have the same energy of the complete graph $K_n.$ A graph $G$ on $n$ vertices is said to be \emph{borderenergetic}  if its energy equals the energy of the complete graph $K_n.$ Some recent results on borderenergetic graphs are the following. 

In \cite{Gutman2015}, it was shown that there exists borderenergetic graphs  on order $n$ for each integer $n\geq 7,$ and all borderenergetic graphs with $7,8,$ and $9$ vertices were determined.
In \cite{JTT2015} it was considered the classes of borderenergetic threshold graphs. For each $n\geq 3,$ it was determined $n-1$ threshold graphs on $n^2$ vertices, pairwise non-cospectral and equienergetic to the complete graph $K_{n^2}.$

Recently,  Hou and Tao  \cite{Hou}, showed that for each $n\geq 2$ and $p\geq 1$  $(p \geq 2$ if $n=2),$ there are  $n-1$ threshold graphs on $pn^2$ vertices, pairwise non-cospectral and equienergetic  with the complete graph $K_{pn^2},$ generalizing the results in \cite{JTT2015}.

In this paper, we continue this investigation in the class of cographs. More precisely, we determine two infinite families of cographs that are borderenergetic.

Here is an outline of  the remainder  of this paper.  In Section 2, we mention the representation of cographs  by a cotree and explain the Diagonalization algorithm. In Section 3,  we determine explicitly the multiplicity $m(\lambda)$ for some classes of cographs,
except $0,-1$ and an upper bound for the remaining eigenvalues.  In Section 4, as application, we present two families of integral non-cospectral and borderenergetic cographs.

\section{Notation and Preliminaries} \label{pre}

Let  $G= (V,E)$ be an undirected graph with vertex set $V$ and edge set $E,$ without loops or multiple edges. We denote the {\em open neighborhood}  of $v,$ by
$$N(v)=\{w|\{v,w\}\in E\}$$ and its {\em closed neighborhood} by $$N[v] = N(v) \cup \{v\}.$$

The {\em adjacency matrix} of $G$, denoted by $A=[a_{ij}]$, is a matrix whose rows and columns are indexed by the vertices of $G$, and is defined to have entries
$$
a_{ij} = \left\{ \begin{array}{rl}
 1 &\mbox{ if $v_{i}v_{j}\in E$} \\
 0 &\mbox{ otherwise}
       \end{array} \right.
$$
A value $\lambda$ is an
{\em eigenvalue} of $G$ if $\det(A - \lambda I_n ) = 0$, and since
$A$ is real and symmetric, its eigenvalues are real numbers.
We denote $m(\lambda)$ the multiplicity of the eigenvalue $\lambda$ of $A.$

\subsection{Cotrees}
A cograph has been rediscovered independently by  several authors
since the 1960's.  Corneil, Lerchs and Burlingham \cite{Stewart} define cographs recursively by the following rules:
\begin{enumerate}
	\item [(i)] a graph on a single vertex is a cograph,
  \item [(ii)] a finite  union of cographs is a cograph,
  \item [(iii)] a finite  join of cographs is a cograph.
\end{enumerate}

In this note, we focus on representing the recursive construction of a cograph using its cotree, that we describe below.

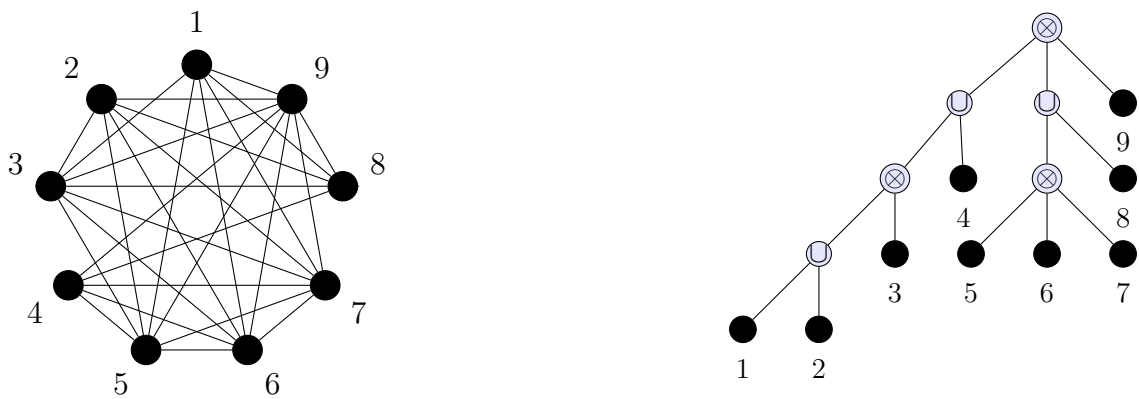
\begin{figure}[h!]
       \begin{minipage}[c]{0.25 \linewidth}
\begin{tikzpicture}
  [scale=0.65,auto=left,every node/.style={circle}]
  \foreach \i/\w in {1/,2/,3/,4/,5/,6/,7/,8/,9/}{
    \node[draw,circle,fill=black,label={360/9 * (\i - 1)+90}:\i] (\i) at ({360/9 * (\i - 1)+90}:3) {\w};} 
  \foreach \from in {3,5,6,7,9}{
    \foreach \to in {1,2,...,\from}
      \draw (\from) -- (\to);}
       \foreach \from in {8}{
       \foreach \to in {1,2,3,4,\from}
       \draw (\from) -- (\to)
       ;}
\end{tikzpicture}
       \end{minipage}\hfill
       \begin{minipage}[l]{0.4 \linewidth}
\begin{tikzpicture}
  [scale=1,auto=left,every node/.style={circle,scale=0.9}]
  \node[draw,circle,fill=black,label=below:$1$] (p) at (-2,3) {};
  \node[draw,circle,fill=black,label=below:$7$] (o) at (3,4) {};
  \node[draw,circle,fill=black,label=below:$6$] (n) at (2,4) {};
   \node[draw,circle,fill=black,label=below:$5$] (q) at (1,4) {};
  \node[draw, circle, fill=blue!10, inner sep=0] (m) at (2,5) {$\otimes$};
  \node[draw,circle,fill=black,label=below:$8$] (l) at (3,5) {};
  \node[draw,circle,fill=black,label=below:$9$] (k) at (3,6) {};
  \node[draw, circle, fill=blue!10, inner sep=0] (j) at (2,6) {$\cup$};
  \node[draw,circle,fill=blue!10, inner sep=0] (h) at (2,7) {$\otimes$};
  \node[draw, circle, fill=blue!10, inner sep=0] (g) at (0.85,6) {$\cup$};
  \node[draw,circle,fill=black,label=below:$4$] (f) at  (0.9,5) {};
  \node[draw,circle,fill=blue!10, inner sep=0] (a) at (0,5) {$\otimes$};
  \node[draw, circle, fill=blue!10, inner sep=0] (b) at (-1,4) {$\cup$};
  \node[draw,circle,fill=black,label=below:$3$] (c) at (0,4) {};
  \node[draw,circle,fill=black,label=below:$2$] (d) at (-1,3) {};
  \path (a) edge node[left]{} (b)
        (a) edge node[below]{} (c)
        (b) edge node[left]{} (d)
        (f) edge node[right]{}(g)
        (g) edge node[left]{}(a)
        (h) edge node[right]{}(j)
        (h) edge node[left]{}(g)
        (h) edge node[left]{}(k)
        (j) edge node[right]{}(l)
        (j) edge node[below]{}(m)
        (m) edge node[below]{}(n)
        (m) edge node[right]{}(o)
        (b) edge node[left]{} (p)
        (m) edge node[left]{} (q);
\end{tikzpicture}
       \end{minipage}
       \caption{The cograph $G=((((v_{1}\cup v_{2})\vee v_{3})\cup v_{4})\vee(((v_{5}\vee v_{6})\vee v_{7})\cup v_{8}))\vee v_{9}$ and its cotree.}
       \label{cotree}
\end{figure}

A cotree $T_G$ of a cograph $G$ is a rooted tree in which any interior vertex $w$ is either of $\cup$ type (corresponding to disjoint union) or $\otimes$ type (corresponding to join). The terminal vertices (leaves) are typeless and represent the vertices of the cograph $G.$  We say that the {\em depth} of the cotree is the number of edges of  the longest path from the root to a leaf. To build a cotree for a connected cograph, we simply place a $\otimes$ at the tree's root, placing $\cup$ on interior vertices with odd depth, and placing $\otimes$ on interior vertices with even depth. To build a cotree for a disconnected cograph, we place $\cup$ at the root, and place $\otimes$ at odd depths, and $\cup$ at even depths. All interior vertices have at least two children. In \cite{BSS2011} this structure is called {\em minimal  cotree}, but  throughout this paper we call it simply a cotree. Figure \ref{cotree} shows  a cograph and its cotree with depth equals  to 4.

Two vertices $u$ and $v$  are {\em duplicates}  if $N(u) = N(v)$ and {\em coduplicates} if $N[u]=N[v].$  In a cograph, any collection of mutually coduplicates (resp. duplicates) vertices, e.g. with the same neighbors and adjacent (resp. not adjacent)  have a common parent of type $ \otimes$ (resp.  $\cup$).
In Figure \ref{cotree}, for example, we have that $v_{1}$ and $v_{2}$ are duplicates because $N(v_{1})=N(v_{2})$, while $v_5, v_6$ and $v_7$ are coduplicates. In fact, a recursive characterization of cographs in terms of the vertex duplication and co-duplication operations is given in \cite{Moha}.

\subsection{Diagonalization}
\label{sub:diag}

An algorithm for constructing a {\em diagonal} matrix congruent to $A + \myvar I$, where $A$ is the adjacency matrix of a cograph,
and $\myvar$ is an arbitrary scalar, using $O(n)$ time and space was developed in \cite{JTT2016}.
This algorithm will be the main tool of this article and, hence, we will make a brief review of the method. For more information, see \cite{JTT2016}.

The algorithm's input  is the cotree $T_G$ and $x.$
Each leaf $v_i, i =1, \ldots,n$ has a value $d_i$ that represents the diagonal element of  $A+xI.$ It initializes  all entries $d_i$ with $x.$ Even though the operations represent rows and columns operations on the matrix $A+xI$, the algorithm is performed on the cotree itself and matrix is never actually used.

In each iteration of the procedure, a pair $\{v_k, v_l\} $ of duplicates or coduplicates vertices with maximum depth  is selected. Then the pair is processed, that is,  assignments are given to $d_k$ and $d_l,$ such that either one or both rows (columns), corresponding to this vertices, are diagonalized. When a $k$ row(column)  corresponding to vertex  $v_k$  has been diagonalized then $v_k$ is  removed from the cotree $T_G,$ it means that $d_k$ has a permanent final value. Then the algorithm moves to the cotree $T_G -v_k.$  The algorithm is shown in Figure~\ref{algo}.

\begin{figure}[h]
{\small
{\tt
\begin{tabbing}
aaa\=aaa\=aaa\=aaa\=aaa\=aaa\=aaa\=aaa\= \kill
     \> INPUT:  cotree $T_G$, scalar $\myvar$\\
     \> OUTPUT: diagonal matrix $D=[d_1, d_2, \ldots, d_n]$ congruent to $A(G) + \myvar I$\\
     \>\\
     \>   $\mbox{ Algorithm}$ Diagonal $(T_{G}, x)$ \\
     \> \> initialize $d_i := \myvar$, for $ 1 \leq i \leq n$ \\
     \> \> {\bf while } $T_G$  has $\geq 2$    leaves      \\
     \> \> \>  select a pair $(v_k, v_l)$  (co)duplicate of maximum depth with  parent $w$\\
     \> \> \>     $\alpha \leftarrow  d_k$    $\beta \leftarrow d_{l}$\\
     \> \> \> {\bf if} $ w=\otimes$\\
     \> \> \> \>  {\bf if} $\alpha + \beta \neq2$  \verb+                //subcase 1a+    \\
     \> \> \> \> \>   $d_{l} \leftarrow \frac{\alpha \beta -1}{\alpha + \beta -2};$ \hspace*{0,25cm} $d_{k} \leftarrow \alpha + \beta -2; $\hspace{0,25cm}   $T_G = T_G - v_k$ \\
     \> \> \> \>  {\bf else if } $\beta=1$ \verb+                //subcase 1b+   \\
     \> \> \> \> \>   $d_{l} \leftarrow 1$ \hspace*{0,25cm}   $d_k  \leftarrow 0;$ \hspace{0,25cm} $T_G = T_G - v_k$ \\
     \> \> \> \>  {\bf else  }  \verb+                      //subcase 1c+   \\
     \> \> \> \> \>   $d_{l} \leftarrow 1$  \hspace*{0,25cm} $d_k \leftarrow -(1-\beta)^2;$ \hspace{0,25cm} $T_G= T_G -v_k;$ \hspace{0,25cm} $T_G = T_G -v_l$  \\
     \> \> \>     {\bf else if} $w=\cup$\\
     \> \> \> \>  {\bf if} $\alpha + \beta \neq 0$  \verb+               //subcase 2a+    \\
     \> \> \> \> \>   $d_{l} \leftarrow \frac{\alpha \beta}{\alpha +\beta};$ \hspace*{0,25cm}   $d_k \leftarrow \alpha +\beta;$ \hspace{0,25cm} $T_G = T_G - v_k$ \\
     \> \> \> \>  {\bf else if } $\beta=0$ \verb+                //subcase 2b+   \\
     \> \> \> \> \>   $d_{l} \leftarrow 0;$ \hspace*{0,25cm}  $d_k  \leftarrow 0;$ \hspace{0,25cm} $T_G = T_G - v_k$ \\
     \> \> \> \>  {\bf else  }  \verb+                      //subcase 2c+   \\
     \> \> \> \> \>   $d_{l} \leftarrow \beta;$  \hspace*{0,25cm} $v_k \leftarrow -\beta;$ \hspace{0,25cm} $T_G =T_G - v_k;$ \hspace{0,25cm} $T_G = T_G - v_l$  \\

     \> \>  {\bf end loop}\\
\end{tabbing}
}}
\caption{\label{algo} Diagonalization algorithm}
\end{figure}

Now, we will present a few results from \cite{JTT2016} that will be used throughout the note. The following theorem is based on  Sylvester's Law of Inertia.

\begin{Thr} \cite{JTT2016}
\label{mainB} Let $D=[d_1,d_2,\ldots,d_n]$ be the diagonal returned  by the diagonalization algorithm $(T_G,-x),$ and assume $D$ has $k_{+}$ positive values, $k_0$ zeros and $k_{-}$ negative  values.
\begin{my_description}
 \item[i]
The number of eigenvalues of $G$ that are greater than $\myvar$ is exactly $k_{+}$.
\item[ii]
The number of eigenvalues of $G$ that are less than $\myvar$ is exactly $k_{-}$.
\item[iii]
The multiplicity of  $\myvar$ is  $k_{0}$.
\end{my_description}
\end{Thr}

The following two lemmas show that, under certain conditions, we can control the assignments made at each iteration.

\begin{Lem}\cite{JTT2016}
\label{lem1}
If $v_1, \ldots, v_m$ have parent $w= \otimes,$ each with diagonal value $y \neq1,$ then the algorithm performs $m-1$ iterations of  \subia~ assigning, during iteration  $j:$
\begin{equation}
d_k  \leftarrow \frac{j+1}{j}(y-1)
\end{equation}
\begin{equation}
d_l  \leftarrow \frac{y+j}{j+1}
\end{equation}
\end{Lem}

\begin{Lem}\cite{JTT2016}
\label{lem2}
If $v_1, \ldots, v_m$ have parent $w= \cup,$ each with diagonal value $y\neq 0,$ then the algorithm performs $m-1$ iterations of  \subiia~ assigning, during iteration  $j:$
\begin{equation}
d_k  \leftarrow \frac{(j+1)}{j}y
\end{equation}
\begin{equation}
d_l  \leftarrow \frac{y}{j+1}
\end{equation}
\end{Lem}

The next three lemmas show that if we start an iteration with some known value then we can control the exit values.
\begin{Lem}\cite{JTT2016}
\label{lem3}
If $\{v_k, v_l \}$ is a pair of coduplicate vertices  processed by Diagonalization with assignments $0 \leq d_k, d_l < 1,$ then $d_k$ becomes permanently negative, and $d_l$  is assigned a value in $(0,1).$
\end{Lem}

\begin{Lem}
\label{lem4}
If $\{v_k, v_l \}$ is a pair of duplicate vertices   processed by Diagonalization with the  assignments $0 < d_k, d_l \leq 1,$ then $d_k$ becomes permanently positive, and $d_l$  is assigned a value in $(0,1).$
\end{Lem}
{\bf Proof:} We notice that the algorithm executes \subiia~, meaning that $d_k=\alpha + \beta >0$ and $d_l= \alpha\beta/(\alpha+\beta).$
The fact that $d_l>0$ is obvious. To see that $d_l < 1$, we observe that if $\alpha = \beta=1$, then $d_k =1/2$. If either (but not both) $\alpha$ or $\beta = 1$, then it is clear that $d_l =\alpha /(\alpha+1) < 1$. Now if $0 < \alpha ,\beta < 1$, then $d_l < 1$ follows from Lemma 3 of \cite{Chang08}.

\begin{Lem}\label{lem5}
During the execution of Diagonalize $(T_G, x)$ with $x \in (0,1),$ all diagonal values of vertices remaining on the cotree are in $(0,1).$
Furthermore, if $d_k$ corresponds to a permanent value of a removed vertex on $T_G -v_k,$ then $d_k \neq 0.$
\end{Lem}
{\bf Proof:}
Let $G$ be a cograph and $T_G$ its cotree.  Initially  all vertices on $T_G$ are in $(0,1).$  Suppose  after  $m$ iterations of Diagonalize all diagonal values of the cotree are in $(0,1)$ and no zero is assigned.  Now consider iteration $m+1$ with a pair $\{ v_k, v_l\}$ and parent $w.$ If $w= \otimes$ then Lemma \ref{lem3}  guarantees  the vertex $d_l$ remaining on the cotree is assigned a value in $(0,1)$ and  the vertex $d_k$ is assigned a permanently negative value.  If $w= \cup$ then Lemma \ref{lem4}  guarantees  the vertex $d_l$ remaining on the cotree is assigned a value in $(0,1)$ and  the vertex $d_k$ is assigned a permanently positive value, completing the proof.

The next result follows from  Lemma \ref{lem5}.

\begin{Thr}  No cograph $G$ has  eigenvalue in the interval $(-1,0).$
\end{Thr}

\begin{figure}[h!]

\begin{center}
\begin{tikzpicture}

  [scale=0.8,auto=left,every node/.style={circle,scale=0.7}]

  \node[draw,circle,fill=black,label=below:$$] (o) at (2.3,4) {};
  \node[draw,circle,fill=black,label=below:$$] (n) at (1.5,4) {};
   \node[draw,circle,fill=black,label=below:$$] (q) at (0.75,4) {};
  \node[draw,circle,fill=black,label=below:$$] (v) at (3,4) {};
  \node[draw,circle,fill=black,label=below:$$] (v2) at (3.5,4) {};
  \node[draw,circle,fill=black,label=below:$$] (z) at (4.2,4) {};
   \node[draw,circle,fill=black,label=below:$$] (z2) at (4.8,4) {};
      \node[draw,circle,fill=black,label=below:$$] (y3) at (5.4,4) {};
         \node[draw,circle,fill=black,label=below:$$] (y4) at (6,4) {};
  \node[draw, circle, fill=blue!10, inner sep=0] (m) at (1.7,5) {$\otimes$};
  \node[draw,circle,fill=blue!10, inner sep=0] (l) at (3,5) {$\otimes$};
  \node[draw,circle,fill=blue!10, inner sep=0] (y) at (4.6,5) {$\otimes$};
   \node[draw,circle,fill=blue!10, inner sep=0] (y2) at (5.4,5) {$\otimes$};
  \node[draw,circle,fill=blue!10, inner sep=0] (k) at (4.6,6) {$\cup$};
  \node[draw, circle, fill=blue!10, inner sep=0] (j) at (2,6) {$\cup$};
  \node[draw,circle,fill=blue!10, inner sep=0] (h) at (2,7) {$\otimes$};
  \node[draw, circle, fill=blue!10, inner sep=0] (g) at (0,6) {$\cup$};
  \node[draw,circle,fill=blue!10, inner sep=0] (f) at  (0.6,5) {$\otimes$};
  \node[draw,circle,fill=blue!10, inner sep=0] (a) at (-1,5) {$\otimes$};
  \node[draw, circle, fill=black,label=below:$$] (b) at (-1.4,4) {};
  \node[draw,circle,fill=black,label=below:$$] (c) at (-0.75,4) {};
  \node[draw,circle,fill=black,label=below:$$] (e) at (0,4) {};

 \path
 (a) edge node[left]{} (b)
        (a) edge node[below]{} (c)

   (f) edge node[below]{} (e)
      (k) edge node[below]{} (y)
       (k) edge node[below]{} (y2)
        (y) edge node[below]{} (z)
         (y) edge node[below]{} (z2)
           (y2) edge node[below]{} (y3)
             (y2) edge node[below]{} (y4)
                   (l) edge node[below]{} (v)
      (l) edge node[below]{} (v2)
        (f) edge node[right]{}(g)
        (g) edge node[left]{}(a)
        (f) edge node[left]{}(q)
        (h) edge node[right]{}(j)
        (h) edge node[left]{}(g)
        (h) edge node[left]{}(k)
        (j) edge node[right]{}(l)
        (j) edge node[below]{}(m)
        (m) edge node[below]{}(n)
        (m) edge node[right]{}(o)
        ;
\end{tikzpicture}
\caption{The cograph with cotree $T_G (3,2,0| 0,0,2).$}
       \label{fig3}
\end{center}
\end{figure}
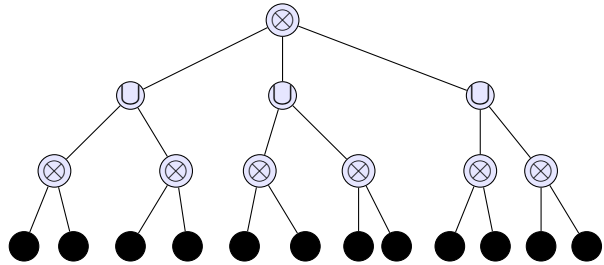

\section{On the multiplicities of eigenvalues in balanced cotrees}
\label{multi}
In this section we study the eigenvalues of cographs that have balanced cotrees.

We say that a cograph $G$  has a \emph{balanced cotree}  $T_G$ with depth $r$  if every interior vertex  with depth $i$ in $T_G$ has the same number of interior vertices and the same number of leaves as direct successors, for $i \in \{1, \ldots,r-1\}$.
We will use the notation $T_G ( a_1,  \ldots, a_r | b_1, \ldots, b_r)$ to represent a balanced cotree of a cograph $G,$
where  the root of $T_G$ has exactly $a_1$ immediate interior vertices and $b_1$ leaves. An interior vertex successor of the root has exactly $a_2$ immediate interior vertices and $b_2$ leaves, and so on.
Thus,  we will assume that $a_1,  \ldots, a_{r-1}$ are positive integers and $a_{r}=0$. Additionally, we assume that $b_1, \ldots, b_{r-1}$ are non negative integer values and $b_{r}\geq 2$. Figure \ref{fig3} shows the balanced cotree $T_G (3,2,0| 0,0,2).$

\subsection{Regular balanced cotrees}
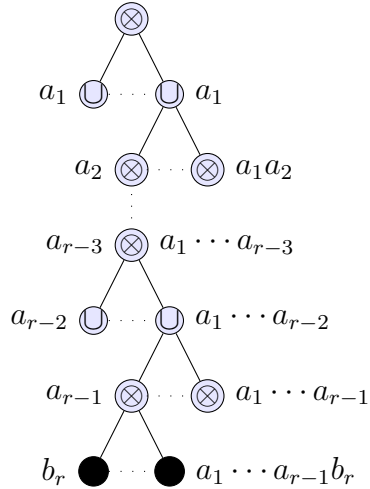
\begin{figure}[h!]
\begin{center}
\begin{tikzpicture}

  [scale=.8,auto=left,every node/.style={circle,scale=0.7}]

\node[draw, circle, fill=blue!10, inner sep=0] (a) at (3,8) {$\otimes$};
\node[draw, circle, fill=blue!10,label=left:$a_{1}$, inner sep=0] (b) at (2.5,7) {$\cup$};
\node[draw, circle, fill=blue!10, label=right:$a_{1}$,inner sep=0] (c) at (3.5,7) {$\cup$};
\draw[loosely dotted] (2.7,7) -- (3.2,7);
\path (a) edge node[left]{} (b);
\path (a) edge node[left]{} (c);

\node[draw, circle, fill=blue!10, label=left:$a_{2}$,inner sep=0] (d) at (3,6) {$\otimes$};
\node[draw, circle, fill=blue!10, label=right:$a_{1}a_{2}$,inner sep=0] (e) at (4,6) {$\otimes$};
\draw[loosely dotted] (3.2,6) -- (3.8,6);
\path (c) edge node[left]{} (d);
\path (c) edge node[left]{} (e);

\draw[loosely dotted] (3,5.2) -- (3,5.8);

\node[draw, circle, fill=blue!10,label=right:$a_{1}\cdots a_{r-3}$,label=left:$a_{r-3}$, inner sep=0] (f) at (3,5) {$\otimes$};
\node[draw, circle, fill=blue!10,label=left:$a_{r-2}$, inner sep=0] (g) at (2.5,4) {$\cup$};
\node[draw, circle, fill=blue!10, label=right:$a_{1}\cdots a_{r-2}$,inner sep=0] (h) at (3.5,4) {$\cup$};
\draw[loosely dotted] (2.7,4) -- (3.2,4);
\path (f) edge node[left]{} (g);
\path (f) edge node[left]{} (h);

\node[draw, circle, fill=blue!10,label=left:$ a_{r-1}$, inner sep=0] (i) at (3,3) {$\otimes$};
\node[draw, circle, fill=blue!10, label=right:$a_{1}\cdots a_{r-1}$,inner sep=0] (j) at (4,3) {$\otimes$};
\draw[loosely dotted] (3.2,3) -- (3.8,3);
\path (h) edge node[left]{} (i);
\path (h) edge node[left]{} (j);

\node[draw,circle,fill=black,label=left:$b_{r}$] (k) at (2.5,2) {};
\node[draw,circle,fill=black,label=right:$a_{1}\cdots a_{r-1}b_{r}$] (l) at (3.5,2) {};

\draw[loosely dotted] (2.7,2) -- (3.3,2);
\path (i) edge node[left]{} (k);
\path (i) edge node[left]{} (l);
\end{tikzpicture}
\caption{Cotree $T_G  (a_1,  \ldots,a_{r-1}, 0| 0, \ldots, 0,b_r)$ with $r$ odd.}
\label{Fig1}
\end{center}
\end{figure}

Here we will study eigenvalues of (regular) cographs $G$ that have balanced cotrees of the type $T_G  (a_1,  \ldots,a_{r-1}, 0| 0, \ldots, 0,b_r)$, whose order is $n =a_{1} a_2  \ldots  a_{r-1} b_r$. We show in Figure \ref{Fig1} a representation of general regular balanced cotree with odd $r$, meaning that level $r-1$ has vertices of type $\otimes$.

The next two theorems are known results and can be found, for example in \cite{BSS2011, JTT2016}.

\begin{Thr}
\label{Tre1}
Let G be a cograph with cotree $T_{G}$ having $\otimes$-nodes $\{w_{1},\ldots,w_{m}\}$, where $w_{i}$ has $t_{i} \geq 1$ terminal children. Then  $m(-1) =\sum_{i=1}^{m}(t_{i}-1)$.
\end{Thr}

\begin{Thr}
\label{Tre2}
Let G be a cograph with cotree $T_{G}$ having $\cup$-nodes $\{w_{1},\ldots,w_{m}\}$, where $w_{i}$ has $t_{i} \geq 1$ terminal children. If $G$ has $j\geq 0$ isolated vertices then  $m(0)= j+\sum_{i=1}^{m}(t_{i}-1)$.
\end{Thr}

Using the above results we can easily prove the next corollary.

\begin{Cor}
\label{Cor1}
Let $G$ be a cograph with balanced cotree $T_G  (a_1,  \ldots,a_{r-1}, 0| 0, \ldots, 0,b_r)$ of order $n =a_{1} a_2  \ldots  a_{r-1} b_r.$
\begin{enumerate}
 \item[(i)]
 If $r$ is odd then $G$ has the eigenvalue $-1$ with multiplicity $a_1 a_2 \ldots a_{r-1}(b_{r} -1)$.
\item[(ii)]
 If $r$ is even then $G$ has the eigenvalue $0$ with multiplicity $a_1 a_2 \ldots a_{r-1}(b_{r} -1)$.
\end{enumerate}
\end{Cor}

\begin{Cor}
\label{Cor2}
Let $G$ be a balanced cotree $T_G  (a_1,  \ldots,a_{r-1}, 0| 0, \ldots,0, b_r)$ of a  cograph $G$ of order $n =a_1 a_2 \ldots a_{r-1} b_r$.
If $r$ is odd (even) then, counting multiplicities, the number of eigenvalues of $G$ other than  $-1$ (0)  is equal to
\begin{equation}
a_1 a_2 \ldots a_{r-1}
\end{equation}
\end{Cor}
 {\bf Proof:} Suppose $r$ is odd. Since $n=a_1 a_2 \ldots a_{r-1} b_r$  and $G$ has $a_1 a_2 \ldots a_{r-1}(b_r -1)$ coduplicates vertices, it  follows that the number of eigenvalues that are distinct from $-1$ is equal to
 $n-a_1 a_2 \ldots a_{r-1}(b_r -1)=a_1 a_2 \ldots a_{r-1}.$ The case $r$ even is similar.

\begin{Lem} \label{lem6} Let $G$ be a cograph with balanced cotree $T_G  (a_1,  \ldots,a_{r-1}, 0| 0, \ldots, 0,b_r)$ of order $n =a_{1} a_2  \ldots  a_{r-1} b_r.$
\begin{enumerate}
 \item[(i)]
   If $r$ is odd then $G$ has the eigenvalue $b_r -1$ with multiplicity $a_1 a_2 \ldots (a_{r-1} -1)$.
\item[(ii)]
If $r$ is even then $G$ has the eigenvalue $-b_r $ with multiplicity $a_1 a_2 \ldots (a_{r-1} -1)$.
\end{enumerate}
\end{Lem}

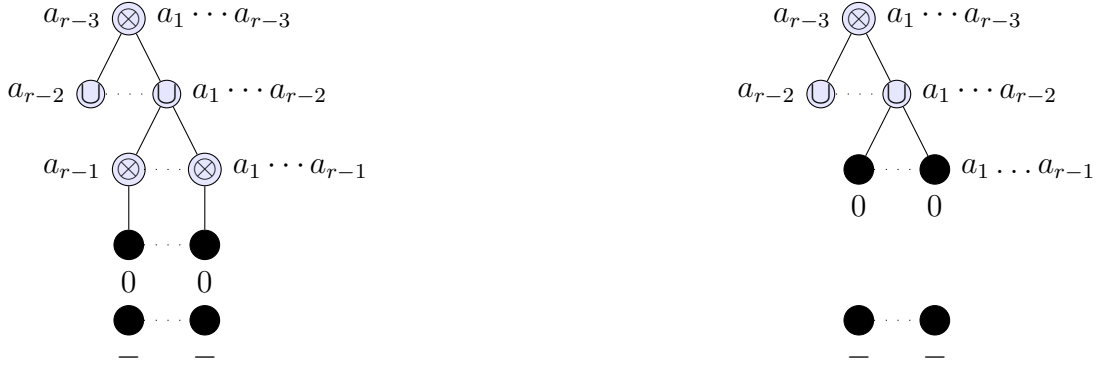
\begin{figure}[t!]
 \begin{minipage}[c]{0.25 \linewidth}
\begin{tikzpicture}

  [scale=,5,auto=left,every node/.style={circle,scale=0.7}]

\node[draw, circle, fill=blue!10,label=right:$a_{1}\cdots a_{r-3}$,label=left:$a_{r-3}$, inner sep=0] (f) at (3,5) {$\otimes$};
\node[draw, circle, fill=blue!10,label=left:$a_{r-2}$, inner sep=0] (g) at (2.5,4) {$\cup$};
\node[draw, circle, fill=blue!10, label=right:$a_{1}\cdots a_{r-2}$,inner sep=0] (h) at (3.5,4) {$\cup$};
\draw[loosely dotted] (2.7,4) -- (3.2,4);
\path (f) edge node[left]{} (g);
\path (f) edge node[left]{} (h);

\node[draw, circle, fill=blue!10,label=left:$ a_{r-1}$, inner sep=0] (i) at (3,3) {$\otimes$};
\node[draw, circle, fill=blue!10, label=right:$a_{1}\cdots a_{r-1}$,inner sep=0] (j) at (4,3) {$\otimes$};
\draw[loosely dotted] (3.2,3) -- (3.8,3);
\path (h) edge node[left]{} (i);
\path (h) edge node[left]{} (j);

\node[draw,circle,fill=black,label=below:$0$] (k) at (3,2) {};
\node[draw,circle,fill=black,label=below:$0$] (l) at (4,2) {};

\path (i) edge node[left]{} (k);
\path (j) edge node[left]{} (l);
\draw[loosely dotted] (3.2,2) -- (3.8,2);

\node[draw,circle,fill=black,label=below:$-$] () at (3,1) {};
\node[draw,circle,fill=black,label=below:$-$] () at (4,1) {};
\draw[loosely dotted] (3.2,1) -- (3.8,1);
\end{tikzpicture}
\end{minipage}\hfill
       \begin{minipage}[l]{0.4 \linewidth}
  \begin{tikzpicture}

  [scale=1,auto=left,every node/.style={circle,scale=0.9}]

\node[draw, circle, fill=blue!10,label=right:$a_{1}\cdots a_{r-3}$,label=left:$a_{r-3}$, inner sep=0] (f) at (3,5) {$\otimes$};
\node[draw, circle, fill=blue!10,label=left:$a_{r-2}$, inner sep=0] (g) at (2.5,4) {$\cup$};
\node[draw, circle, fill=blue!10, label=right:$a_{1}\cdots a_{r-2}$,inner sep=0] (h) at (3.5,4) {$\cup$};
\draw[loosely dotted] (2.7,4) -- (3.2,4);
\path (f) edge node[left]{} (g);
\path (f) edge node[left]{} (h);

\node[draw,circle,fill=black,label=below:$0$] (i) at (3,3) {};
\node[draw,circle,fill=black,label=below:$0$,label=right:$a_{1}\ldots a_{r-1}$] (j) at (4,3) {};
\draw[loosely dotted] (3.2,3) -- (3.8,3);
\path (h) edge node[left]{} (i);
\path (h) edge node[left]{} (j);

\node[draw,circle,fill=black,label=below:$-$] () at (3,1) {};
\node[draw,circle,fill=black,label=below:$-$] () at (4,1) {};
\draw[loosely dotted] (3.2,1) -- (3.8,1);

\end{tikzpicture}
\end{minipage}
\caption{Processing deepest  $\bigotimes$ level.}
       \label{figure11}
\end{figure}

\noindent{\bf Proof:}
We assume that $r$ is odd. The case even is similar.
Consider  $ x= -(b_r -1)$ and execute the algorithm Diagonalization with input $(T_G, x).$
By Theorem \ref{mainB} we have to prove that the algorithm creates $a_1 a_2 \ldots (a_{r-1} -1)$ null permanent values.
Since $G$ has coduplicate vertices, see Figure \ref{Fig1}, and $ x\neq 1,$ we apply Lemma \ref{lem1} and  after $b_r -1$ iterations for each $\otimes$ vertice at level $r-1$ , the remaining vertices on the cotree receive
$$ d_l  \leftarrow \frac{ -(b_r -1) + b_r -1}{b_r -1 +1} =0,$$
and the removed vertices receive
$$d_k \leftarrow -\frac{j+1}{j} b_r<0,\mbox{ for } j = 1, \ldots, b_r -1.$$
This is illustrated on the left of Figure \ref{figure11}.

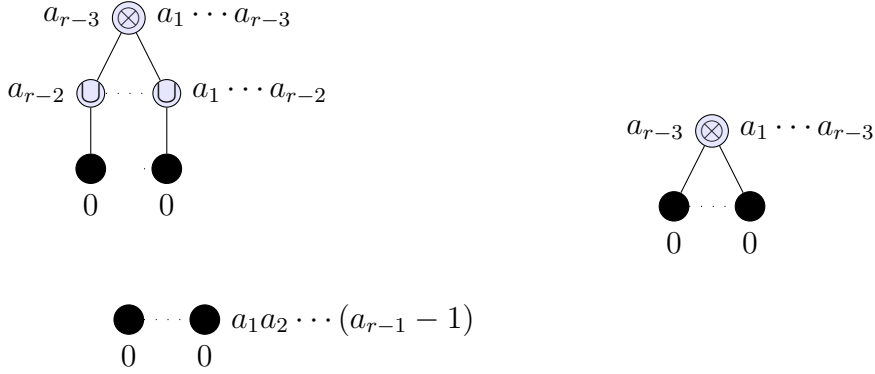
\begin{figure}[h!]
\begin{minipage}[r]{0.5 \linewidth}
\begin{tikzpicture}

  [scale=1,auto=left,every node/.style={circle,scale=0.9}]

\node[draw, circle, fill=blue!10,label=right:$a_{1}\cdots a_{r-3}$,label=left:$a_{r-3}$, inner sep=0] (f) at (3,5) {$\otimes$};
\node[draw, circle, fill=blue!10,label=left:$a_{r-2}$, inner sep=0] (g) at (2.5,4) {$\cup$};
\node[draw, circle, fill=blue!10, label=right:$a_{1}\cdots a_{r-2}$,inner sep=0] (h) at (3.5,4) {$\cup$};
\draw[loosely dotted] (2.7,4) -- (3.2,4);
\path (f) edge node[left]{} (g);
\path (f) edge node[left]{} (h);

\node[draw,circle,fill=black,label=below:$0$] (i) at (2.5,3) {};
\node[draw,circle,fill=black,label=below:$0$] (j) at (3.5,3) {};
\draw[loosely dotted] (3.2,3) -- (3.8,3);
\path (g) edge node[left]{} (i);
\path (h) edge node[left]{} (j);

\node[draw,circle,fill=black,label=below:$0$] () at (3,1) {};
\node[draw,circle,fill=black,label=below:$0$,label=right:$a_{1}a_{2}\cdots (a_{r-1}-1)$] () at (4,1) {};
\draw[loosely dotted] (3.2,1) -- (3.8,1);
\end{tikzpicture}
\end{minipage}
\begin{minipage}[l]{0.2 \linewidth}

\begin{tikzpicture}

  [scale=1,auto=left,every node/.style={circle,scale=0.9}]

\node[draw, circle, fill=blue!10,label=right:$a_{1}\cdots a_{r-3}$,label=left:$a_{r-3}$, inner sep=0] (f) at (3,5) {$\otimes$};
\node[draw,circle,fill=black,label=below:$0$] (g) at (2.5,4) {};
\node[draw,circle,fill=black,label=below:$0$] (h) at (3.5,4) {};
\draw[loosely dotted] (2.7,4) -- (3.2,4);
\path (f) edge node[left]{} (g);
\path (f) edge node[left]{} (h);
\end{tikzpicture}
\end{minipage}
\caption{Processing the deepest $\bigcup$ level}
       \label{figure12}
\end{figure}

Now the leaves at level $r$ move up to the $\cup$ vertices as on the right of Figure \ref{figure11} and we process them. Notice that we have duplicate leaves with null value. Then the algorithm performs \subiib~ at the leaves in each vertex $\cup$ and it creates $a_1 a_2 \ldots (a_{r-1}-1)$ permanent zeros in the removed vertices. The remaining vertices keep the value zero, as shown on the left of Figure \ref{figure12}. So $m(b_r -1) \geq a_1 a_2 \ldots (a_{r-1} -1).$

Now we show that no more permanent zeros are created. The zero value vertices now move up to the next $\bigotimes$ level.
Notice that, see right of Figure \ref{figure12}, we have coduplicate vertices in the remaining tree with assignments equal to $0$.
Using Lemma \ref{lem3} once and then Lemma \ref{lem5}, we know that no null value will be generated and it proves that $m(b_r -1) = a_1 a_2 \ldots (a_{r-1} -1).$

In the next theorem we present a bound for the eigenvalues of regular balanced cotrees.

\begin{Thr}
\label{main3}
Let $G$ be a cograph with balanced cotree $T_G  (a_1,  \ldots,a_{r-1}, 0| 0, \ldots, 0,b_r)$ order $n =a_{1} a_2  \ldots  a_{r-1} b_r.$
\begin{enumerate}
	\item [(i)] If $r$ is odd and $\lambda\neq -1$, $b_{r}-1$ then $m(\lambda)\leq a_{1}\cdots a_{r-2}$;
	\item [(ii)] If $r$ is even and $\lambda\neq 0$, $-b_{r}$ then $m(\lambda)\leq a_{1}\cdots a_{r-2}$.
\end{enumerate}
\end{Thr}
{\bf Proof:} Suppose that $r$ is odd. Then $m(\lambda) \leq n - m(-1) - m(b_r-1)=a_1 a_2 \ldots a_{r-1} b_r - a_1 a_2 \ldots a_{r-1}(b_{r} -1) - a_1 a_2 \ldots (a_{r-1} -1)= a_{1}\cdots a_{r-2}$.

\subsection{Non-regular balanced cotrees}

We now define two types of cotrees depending on whether its depth $r$ is even or odd.
Let  $T_G (a_1, \ldots, a_{r-1},0| b_1, b_2,\ldots, b_r)$  be a balanced cotree  defined as follows:

If $r$ is even then, for $1 \leq i \leq r-1$,
$\left\{
\begin{array}{lr}
b_{i} =0 & \mbox{if $i$ is odd };\\
b_{i}\geq b_{r} & \mbox{if $i$ is  even.}
\end{array} \right.$

If $r$ is odd then, for $1 \leq i \leq r-1$,
$\left\{
\begin{array}{lr}
b_{i} =0 &\mbox{if $i$ is even}; \\
b_{i}\geq b_{r}  &\mbox{if $i$ is  odd.}
\end{array} \right.$

\begin{figure}[h!]
       \begin{minipage}[c]{0.25 \linewidth}
\begin{tikzpicture}
  [scale=1,auto=left,every node/.style={circle,scale=0.9}]
\node[draw, circle, fill=blue!10,label=left:$a_{r-3}$, inner sep=0] (1) at (3,3) {$\cup$};
\node[draw, circle, fill=blue!10,label=left:$a_{r-2}$, inner sep=0] (2) at (2,2) {$\otimes$};
\node[draw, circle, fill=blue!10, inner sep=0] (3) at (3,2) {$\otimes$};
\draw[loosely dotted] (2.2,2) -- (2.8,2);
\node[draw,circle,fill=black] (4) at (3.5,2) {};
\node[draw,circle,fill=black,label=right:$b_{r-2}$] (5) at (4.5,2) {};
\draw[loosely dotted] (3.7,2) -- (4.3,2);

\node[draw, circle, fill=blue!10,label=left:$a_{r-1}$, inner sep=0] (6) at (2.5,1) {$\cup$};
\node[draw, circle, fill=blue!10, inner sep=0] (7) at (3.5,1) {$\cup$};
\draw[loosely dotted] (2.7,1) -- (3.3,1);

\node[draw,circle,fill=black] (8) at (2,0) {};
\node[draw,circle,fill=black] (9) at (2.8,0) {};
\draw[loosely dotted] (2.2,0) -- (2.6,0);

\node[draw,circle,fill=black] (10) at (3.2,0) {};
\node[draw,circle,fill=black,label=right:$b_{r}$] (11) at (4,0) {};
\draw[loosely dotted] (3.4,0) -- (3.8,0);

\path (1) edge node[left]{} (2);
\path (1) edge node[left]{} (3);
\path (1) edge node[left]{} (4);
\path (1) edge node[left]{} (5);
\path (3) edge node[left]{} (6);
\path (3) edge node[left]{} (7);
\path (6) edge node[left]{} (8);
\path (6) edge node[left]{} (9);
\path (7) edge node[left]{} (10);
\path (7) edge node[left]{} (11);
\end{tikzpicture}
\end{minipage}\hfil
\begin{minipage}[r]{.4 \linewidth}
  \begin{tikzpicture}
  [scale=1,auto=left,every node/.style={circle,scale=0.9}]
\node[draw, circle, fill=blue!10,label=left:$a_{r-3}$, inner sep=0] (1) at (3,3) {$\cup$};
\node[draw, circle, fill=blue!10,label=left:$a_{r-2}$, inner sep=0] (2) at (2,2) {$\otimes$};
\node[draw, circle, fill=blue!10, inner sep=0] (3) at (3,2) {$\otimes$};
\draw[loosely dotted] (2.2,2) -- (2.8,2);
\node[draw,circle,fill=black] (4) at (3.5,2) {};
\node[draw,circle,fill=black,label=right:$b_{r-2}$] (5) at (4.5,2) {};
\draw[loosely dotted] (3.7,2) -- (4.3,2);

\node[draw, circle, fill=blue!10,label=left:$a_{r-1}$, inner sep=0] (6) at (2.5,1) {$\cup$};
\node[draw, circle, fill=blue!10, inner sep=0] (7) at (3.5,1) {$\cup$};
\draw[loosely dotted] (2.7,1) -- (3.3,1);

\node[draw,circle,fill=black,label=below:$1$] (8) at (2.5,0) {};

\node[draw,circle,fill=black,label=below:$1$] (10) at (3.5,0) {};

\path (1) edge node[left]{} (2);
\path (1) edge node[left]{} (3);
\path (1) edge node[left]{} (4);
\path (1) edge node[left]{} (5);
\path (3) edge node[left]{} (6);
\path (3) edge node[left]{} (7);
\path (6) edge node[left]{} (8);
\path (7) edge node[left]{} (10);
\end{tikzpicture}
\end{minipage}
\caption{Processing the deepest $\bigcup$ level}
       \label{figure19}
\end{figure}
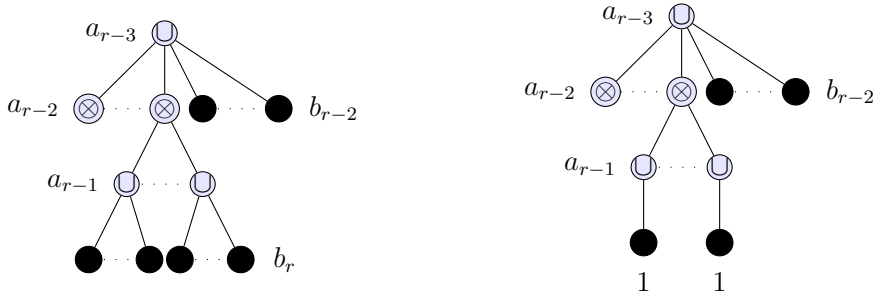

\begin{Thr}
\label{lem7}
Let $G$ be a cograph with balanced cotree $T_G  (a_1,  \ldots,a_{r-1}, 0| b_1, b_2, \ldots,b_r)$ defined above.
\begin{enumerate}
 \item[(i)]
 If $r$ is odd then $G$ has the eigenvalue $b_r -1$ with multiplicity $a_1 a_2 \ldots (a_{r-1} -1)$.
\item[(ii)]
If $r$ is even then $G$ has the eigenvalue $-b_r $ with multiplicity $a_1 a_2 \ldots (a_{r-1} -1)$.
\end{enumerate}
\end{Thr}
{\bf Proof:}
We assume that $r$ is even. The case odd is similar. The illustration of the initial configuration is given on the left of Figure \ref{figure19}. Consider  $ x= -b_r $ and execute the algorithm Diagonalization with input $(T_G, x).$ By Theorem \ref{mainB} we have to prove that the algorithm creates at least $a_1 a_2 \ldots (a_{r-1} -1)$ permanent null values. Applying Lemma \ref{lem2} at each vertex $\cup$ at level $r-1$, the following assignments are made
$$
\begin{array}{cccc}
  d_{k} & \leftarrow & \frac{(j+1)}{j}b_{r}>0, & j=1,\ldots,b_{r}-1; \\
  d_{l} & \leftarrow & \frac{b_{r}}{b_{r}-1+1}=1. &
\end{array}
$$
The removed leaves have a permanent positive value and the remaining vertices have value $1$, as illustrated on the right of Figure \ref{figure19}.

\

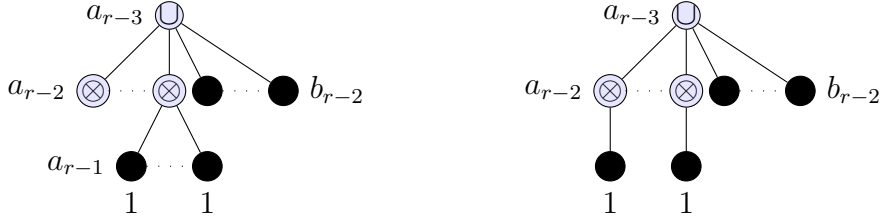
\begin{figure}[h!]
\begin{minipage}[c]{.25 \linewidth}
\begin{tikzpicture}

 [scale=1,auto=left,every node/.style={circle,scale=0.9}]

\node[draw, circle, fill=blue!10,label=left:$a_{r-3}$, inner sep=0] (1) at (3,3) {$\cup$};
\node[draw, circle, fill=blue!10,label=left:$a_{r-2}$, inner sep=0] (2) at (2,2) {$\otimes$};
\node[draw, circle, fill=blue!10, inner sep=0] (3) at (3,2) {$\otimes$};
\draw[loosely dotted] (2.2,2) -- (2.8,2);
\node[draw,circle,fill=black] (4) at (3.5,2) {};
\node[draw,circle,fill=black,label=right:$b_{r-2}$] (5) at (4.5,2) {};
\draw[loosely dotted] (3.7,2) -- (4.3,2);

\node[draw, circle, fill=black,label=left:$a_{r-1}$,label=below:$1$] (6) at (2.5,1) {};
\node[draw, circle, fill=black,label=below:$1$] (7) at (3.5,1) {};
\draw[loosely dotted] (2.7,1) -- (3.3,1);

\path (1) edge node[left]{} (2);
\path (1) edge node[left]{} (3);
\path (1) edge node[left]{} (4);
\path (1) edge node[left]{} (5);
\path (3) edge node[left]{} (6);
\path (3) edge node[left]{} (7);
\end{tikzpicture}
\end{minipage}\hfil
\begin{minipage}[r]{.4 \linewidth}
\begin{tikzpicture}

  [scale=1,auto=left,every node/.style={circle,scale=0.9}]s

\node[draw, circle, fill=blue!10,label=left:$a_{r-3}$, inner sep=0] (1) at (3,3) {$\cup$};
\node[draw, circle, fill=blue!10,label=left:$a_{r-2}$, inner sep=0] (2) at (2,2) {$\otimes$};
\node[draw, circle, fill=blue!10, inner sep=0] (3) at (3,2) {$\otimes$};
\draw[loosely dotted] (2.2,2) -- (2.8,2);
\node[draw,circle,fill=black] (4) at (3.5,2) {};
\node[draw,circle,fill=black,label=right:$b_{r-2}$] (5) at (4.5,2) {};
\draw[loosely dotted] (3.7,2) -- (4.3,2);
\node[draw, circle, fill=black,label=below:$1$] (6) at (2,1) {};
\node[draw, circle, fill=black,label=below:$1$] (7) at (3,1) {};

\path (1) edge node[left]{} (2);
\path (1) edge node[left]{} (3);
\path (1) edge node[left]{} (4);
\path (1) edge node[left]{} (5);
\path (2) edge node[left]{} (6);
\path (3) edge node[left]{} (7);
\end{tikzpicture}
\end{minipage}
\caption{Processing the deepest $\bigotimes$ level.}
       \label{figure21}
\end{figure}

Now the vertices remaining (with value 1) are moved up and become leaves of a $\bigotimes$ vertex, as seeing on the left of Figure \ref{figure21}. We perform subcase 1b and then the $a_{1}\cdots a_{r-2}(a_{r-1}-1)$ removed vertices receive the value $0$ and the remaining leaves receive $1$ as shown on the right of Figure \ref{figure21}, and so $m(-b_{r})\geq a_1 a_2 \ldots (a_{r-1} -1)$.

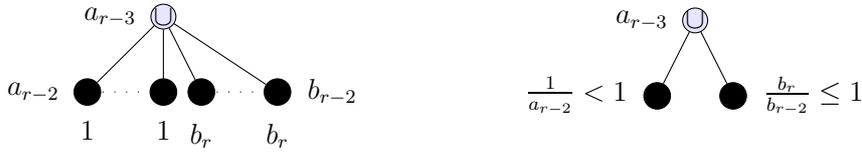
\begin{figure}[h!]
\begin{minipage}[c]{.25 \linewidth}

\begin{tikzpicture}
  [scale=1,auto=left,every node/.style={circle,scale=0.9}]
\node[draw, circle, fill=blue!10,label=left:$a_{r-3}$, inner sep=0] (1) at (3,3) {$\cup$};
\node[draw, circle, fill=black,label=left:$a_{r-2}$,label=below:$1$] (2) at (2,2) {};
\node[draw, circle, fill=black,label=below:$1$] (3) at (3,2) {};
\draw[loosely dotted] (2.2,2) -- (2.8,2);
\node[draw,circle,fill=black,label=below:$b_{r}$] (4) at (3.5,2) {};
\node[draw,circle,fill=black,label=right:$b_{r-2}$,label=below:$b_{r}$] (5) at (4.5,2) {};
\draw[loosely dotted] (3.7,2) -- (4.3,2);
\path (1) edge node[left]{} (2);
\path (1) edge node[left]{} (3);
\path (1) edge node[left]{} (4);
\path (1) edge node[left]{} (5);
\end{tikzpicture}
\end{minipage}\hfil
\begin{minipage}[r]{.4 \linewidth}

\begin{tikzpicture}
  [scale=1,auto=left,every node/.style={circle,scale=0.9}]
\node[draw, circle, fill=blue!10,label=left:$a_{r-3}$, inner sep=0] (1) at (3,3) {$\cup$};
\node[draw, circle, fill=black,label=left:$\frac{1}{a_{r-2}}<1$] (2) at (2.5,2) {};
\node[draw, circle, fill=black,label=right:$\frac{b_{r}}{b_{r-2}}\leq1$] (3) at (3.5,2) {};
\path (1) edge node[left]{} (2);
\path (1) edge node[left]{} (3);
\end{tikzpicture}
\end{minipage}
\caption{Processing level $r-2$.}
       \label{figure23}
\end{figure}

Now the vertices with value $1$ move to level $r-2$ as shown on the left of Figure \ref{figure23}. At each vertex $\cup$ at level $r-3$ we start processing the vertices with value $1$, and by Lemma \ref{lem2}:
$$
\begin{array}{cccc}
  d_{k} & \leftarrow & \frac{(j+1)}{j}1>0, & j=1,\ldots,a_{r-2}-1; \\
  d_{l} & \leftarrow & \frac{1}{a_{r-2}-1+1}=\frac{1}{a_{r-2}}. &
\end{array}
$$
Then we process the vertices with value $b_{r}$ using Lemma \ref{lem2}:
$$
\begin{array}{cccc}
  d_{k} & \leftarrow & \frac{(j+1)}{j}b_{r}>0, & j=1,\ldots,b_{r-2}-1; \\
  d_{l} & \leftarrow & \frac{b_{r}}{b_{r-2}-1+1}=\frac{b_{r}}{b_{r-2}} &
\end{array}
$$
The right of Figure \ref{figure23} represents the last iteration in each vertex $\cup$ at level $r-2$. Notice that each remaining leaf has a value in $(0,1]$ and using the same argument as in Lemma \ref{lem6}, we can prove that  no more zeros are assigned and the remaining vertices on the cotree are in $(0,1)$, proving that $m(-b_r)=a_1 a_2 \ldots (a_{r-1} -1)$.

\section{Borderenergetic Cographs}

In this section we present some families of non-cospectral and borderenergetic cographs.

\begin{figure}[h!]
\begin{center}
\begin{tikzpicture}
  [scale=1,auto=left,every node/.style={circle,scale=0.9}]
  \node[draw,circle,fill=blue!10, inner sep=0] (a) at (2,2) {$\otimes$};
  \node[draw, circle, fill=black,label=left:$\ldots$, label=right:$a $] (b) at (2,1) {};
  \node[draw, circle, fill=blue!10, inner sep=0] (f) at (3.5,1) {$\cup$};
  \node[draw,circle,fill=black, label=left:$1$] (c) at (0.75,1) {};
  \node[draw,circle,fill=black,  label=below:$b$] (d) at (3,-1) {};
  \node[draw,circle,fill=black,  label=below:$1$,  label=right:$\ldots$] (e) at (2,-1) {};
  \node[draw,circle,fill=black, label=below:$1$,  label=right:$\ldots$] (i) at (4,-1) {};
  \node[draw,circle,fill=black,  label=below:$b$] (j) at (5,-1) {};

  \node[draw,circle,fill=blue!10, inner sep=0, label=left:$1$, label=right:$\ldots$] (g) at (2.5,0) {$\otimes$};
  \node[draw,circle,fill=blue!10, inner sep=0,label=left:$\ldots$, label=right:$(a-1)(b-1)$] (h) at (4.5,0) {$\otimes$};
  \path (a) edge node[left]{} (b)
        (a) edge node[below]{} (c)
        (a) edge node[below]{} (f)
      (h) edge node[below]{} (i)
      (h) edge node[below]{}  (j)
       (g) edge node[below]{} (d)
      (g) edge node[below]{} (e)
        (f) edge node[left]{} (g)
        (f) edge node[left]{} (h);
\end{tikzpicture}
       \caption{The cotree $T_{G}$}
       \label{Figu5}
\end{center}
\end{figure}
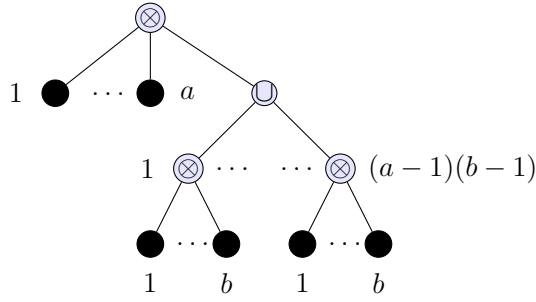

Consider the cograph $G=K_{a} \otimes (a-1)(b-1)K_{b}$, of order $n=a+b(a-1)(b-1)$. We observe that $G$ has the balanced cotree $T_{G}(1, (a-1)(b-1), 0| a,0,b)$, represented in Figure \ref{Figu5}.

\begin{Lem}
\label{spec1}
Let $G=K_{a} \otimes (a-1)(b-1)K_{b}$ be the cograph $G$ of Figure \ref{Figu5} of order $n=a+b(a-1)(b-1)$, for fixed values $a\geq b\geq 2$. The spectrum of $G$ is
$$-(a-1)(b-1); -1; b-1; ab-1$$ with multiplicity $$1; (a-1)[(b-1)^{2}+1]; (a-1)(b-1)-1; 1,$$ respectively.
\end{Lem}
{\bf Proof:}
Using Theorem \ref{Tre1} with $m=(a-1)(b-1)+1$, $t_{1}=\cdots=t_{m-1}=b$ and $t_{m}=a$. We compute the multiplicity of $-1$:
$$m(-1)=\sum_{i=1}^{m}(t_{i}-1)=(m-1)(b-1)+(a-1)=(a-1)[(b-1)^{2}+1].$$ Since $T_G(1,(a-1)(b-1),0|a,0,b)$ has a non regular balanced cotree, by Lemma \ref{lem7}, $b-1$ is an eigenvalue with multiplicity $(a-1)(b-1)-1.$
Now, we will prove, by Theorem \ref{mainB}, that $ab-1$ is an eigenvalue of multiplicity $1$ by showing that the algorithm D Diagonalize with input $(T_{G},-ab+1)$, creates a single zero in the $T_G$.

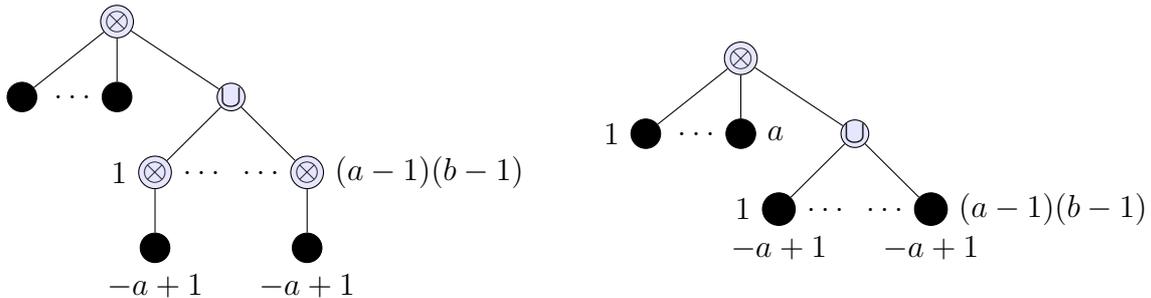
\begin{figure}[h!]
 \begin{minipage}[l]{.4 \linewidth}

\begin{tikzpicture}

  [scale=.9,auto=left,every node/.style={circle,scale=0.9}]

  \node[draw,circle,fill=blue!10, inner sep=0] (a) at (2,2) {$\otimes$};
  \node[draw, circle, fill=black,label=left:$\ldots$, label=right:$ $] (b) at (2,1) {};
  \node[draw, circle, fill=blue!10, inner sep=0] (f) at (3.5,1) {$\cup$};
  \node[draw,circle,fill=black, label=left:$$] (c) at (0.75,1) {};

  \node[draw,circle,fill=black,  label=below:$-a+1$] (d) at (2.5,-1) {};
 (i) at (4,-1) {};
  \node[draw,circle,fill=black,  label=below:$-a+1$] (j) at (4.5,-1) {};
  \node[draw,circle,fill=blue!10, inner sep=0, label=left:$1$, label=right:$\ldots$] (g) at (2.5,0) {$\otimes$};
  \node[draw,circle,fill=blue!10, inner sep=0,label=left:$\ldots$, label=right:$(a-1)(b-1)$] (h) at (4.5,0) {$\otimes$};
  \path (a) edge node[left]{} (b)
        (a) edge node[below]{} (c)
        (a) edge node[below]{} (f)
      (h) edge node[below]{}  (j)
       (g) edge node[below]{} (d)
        (f) edge node[left]{} (g)
        (f) edge node[left]{} (h);
\end{tikzpicture}
  \end{minipage}\hfil
\begin{minipage}[r]{.4 \linewidth}
\begin{tikzpicture}

  [scale=.9,auto=left,every node/.style={circle,scale=0.9}]

  \node[draw,circle,fill=blue!10, inner sep=0] (a) at (2,2) {$\otimes$};
  \node[draw, circle, fill=black,label=left:$\ldots$, label=right:$a $] (b) at (2,1) {};
  \node[draw, circle, fill=blue!10, inner sep=0] (f) at (3.5,1) {$\cup$};
  \node[draw,circle,fill=black, label=left:$1$] (c) at (0.75,1) {};

  \node[draw,circle,fill=black, inner sep=0,label=below:$-a+1$, label=left:$1$, label=right:$\ldots$] (g) at (2.5,0) {$\otimes$};
  \node[draw,circle,fill=black, inner sep=0,label=below:$-a+1$,label=left:$\ldots$, label=right:$(a-1)(b-1)$] (h) at (4.5,0) {$\otimes$};
  \path (a) edge node[left]{} (b)
        (a) edge node[below]{} (c)
        (a) edge node[below]{} (f)
        (f) edge node[left]{} (g)
        (f) edge node[left]{} (h);
\end{tikzpicture}
\end{minipage}
 \caption{Processing deepest level}
 \label{figure30}
\end{figure}

We initialize the leaves with value $-ab+1\neq 1$. Then we can use Lemma \ref{lem1}, and for each  $\otimes$ vertex, we have that
$$
\begin{array}{cccc}
  d_{k} & \leftarrow & \frac{(j+1)}{j}(-ab), & j=1,\ldots,b-1; \\
  d_{l} & \leftarrow & -a+1&
\end{array}
$$
where $d_{k}$ represents the removed leaves and $d_{l}$ the remaining ones. The left of Figure \ref{figure30} represents the cotree yet to be processed. Now, the leaves at depth 3 move up to the $\cup$ vertices at depth 2, as on the right of Figure \ref{figure30}.

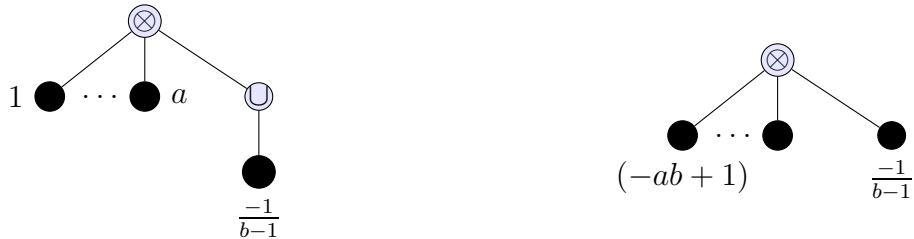
\begin{figure}[h!]
\begin{minipage}[r]{.4 \linewidth}
\begin{tikzpicture}

  [scale=1,auto=left,every node/.style={circle,scale=0.9}]

  \node[draw,circle,fill=blue!10, inner sep=0] (a) at (2,2) {$\otimes$};
  \node[draw, circle, fill=black,label=left:$\ldots$, label=right:$a $] (b) at (2,1) {};
  \node[draw, circle, fill=blue!10, inner sep=0] (f) at (3.5,1) {$\cup$};
  \node[draw,circle,fill=black, label=left:$1$] (c) at (0.75,1) {};

  \node[draw,circle,fill=black, inner sep=0,label=below:$\frac{-1}{b-1}$, label=left:$$, label=right:$$] (g) at (3.5,0) {$\otimes$};

  \path (a) edge node[left]{} (b)
        (a) edge node[below]{} (c)
        (a) edge node[below]{} (f)

        (f) edge node[left]{} (g);
\end{tikzpicture}
\end{minipage}\hfil
\begin{minipage}[r]{.4 \linewidth}
\begin{tikzpicture}

  [scale=1,auto=left,every node/.style={circle,scale=0.9}]

  \node[draw,circle,fill=blue!10, inner sep=0] (a) at (2,2) {$\otimes$};
  \node[draw, circle, fill=black,label=left:$\ldots$] (b) at (2,1) {};
  \node[draw, circle, fill=black,label=below:$\frac{-1}{b-1}$, inner sep=0] (f) at (3.5,1) {$\cup$};
  \node[draw,circle,fill=black, label=below:$(-ab+1)$] (c) at (0.75,1) {};

  \path (a) edge node[left]{} (b)
        (a) edge node[below]{} (c)
        (a) edge node[below]{} (f);
\end{tikzpicture}
\end{minipage}
       \caption{Processing level 2}
       \label{figure32}
\end{figure}

In the next step we use Lemma \ref{lem2} because the duplicate vertices at depth 2 have assignments equal to $-a+1\neq0$. We obtain
$$
\begin{array}{cccc}
  d_{k} & \leftarrow & \frac{(j+1)}{j}(-a+1), & j=1,\ldots,(a-1)(b-1)-1; \\
  d_{l} & \leftarrow & \frac{-1}{b-1}.&
\end{array}
$$
As the left of Figure \ref{figure32} shows, the remaining leaf at depth 2 moves up to depth 1, as on the right of Figure \ref{figure32}.

At depth one, there are  $a+1$ coduplicate vertices. $a$ with value $-ab+1$ and one with value $\frac{ -1}{b-1}$ as the right of Figure \ref{figure32}. The algorithm processes, by Lemma \ref{lem1}, the leaves with value $-ab+1$ first and it generates the following assignments
$$
\begin{array}{cccc}
  d_{k} & \leftarrow & \frac{(j+1)}{j}(-ab), & j=1,\ldots,(a-1); \\
  d_{l} & \leftarrow & -b+1.&
\end{array}
$$

 The last step of the algorithm is to process the two remaining vertices whose values are $\alpha = -b+1$ and $\beta = \frac{-1}{b-1}.$ Since $\alpha, \beta<0$ then the algorithm performs \subia~ and assigns

 $$d_k \leftarrow \alpha + \beta -2 = \frac{-b^{2}}{b-1}\mbox{ and } d_l \leftarrow \frac{\alpha \beta -1}{\alpha +\beta -2} =0,$$
creating a negative value and a zero for the last two diagonal entries, so $m(ab-1)=1$.

Using the fact that sum of eigenvalues must be zero, we obtain the remaining  eigenvalue $-(a-1)(b-1)$ of $G$, proving the result.

The following theorem follows directly from Lemma \ref{spec1} and summarizes the results for the family of cographs represented in Figure \ref{Figu5}.

\begin{Thr}
Let $G=K_{a} 	\otimes (a-1)(b-1)K_{b}$ be the cograph of order $n=a+b(a-1)(b-1)$, for fixed values $a\geq b\geq 2$ represented in Figure \ref{Figu5}. Then $G$ is an integral cograph, non-cospectral and borderenergetic to $K_{n}$.
\end{Thr}
{\bf Proof:}
It is well known that the Spec$(K_{n})=\{(-1)^{n-1},(n-1)^{1}\}$ and and, hence, $E(K_n) = 2(n-1)$. Using Lemma \ref{spec1} we can compute the energy of $G$ as follows
$$E(G)=(a-1)(b-1)+(1)(a-1)[(b-1)^{2}+1]+(b-1)[(a-1)(b-1)-1]+(ab-1)=2(n-1).$$

\begin{figure}[h!]
\begin{center}
\begin{tikzpicture}

  [scale=1,auto=left,every node/.style={circle,scale=0.9}]

  \node[draw,circle,fill=blue!10, inner sep=0] (a) at (2,2) {$\otimes$};
  \node[draw, circle, fill=blue!10, inner sep=0] (f) at (4,1) {$\cup$};
  \node[draw,circle,fill=blue!10, inner sep=0] (c) at (0.5,1) {$\cup$};
  \node[draw,circle,fill=black] (d) at (4,-1) {};
  \node[draw,circle,fill=black] (e) at (3,-1) {};
  \node[draw,circle,fill=black] (i) at (4.5,-1) {};
  \node[draw,circle,fill=black] (j) at (5.5,-1) {};

  \node[draw,circle,fill=black] (m) at (-0.5,-1) {};
   \node[draw,circle,fill=black] (n) at (0.5,-1) {};
    \node[draw,circle,fill=black] (o) at (1,-1) {};
   \node[draw,circle,fill=black] (p) at (2,-1) {};

  \node[draw,circle,fill=blue!10, inner sep=0, label=left:$1$, label=right:$\ldots$] (k) at (0,0) {$\otimes$};
  \node[draw,circle,fill=blue!10, inner sep=0, label=left:$\ldots$, label=right:$p+1$] (b) at (1.5,0) {$\otimes$};
  \node[draw,circle,fill=blue!10, inner sep=0, label=left:$1$, label=right:$\ldots$] (g) at (3.5,0) {$\otimes$};
  \node[draw,circle,fill=blue!10, inner sep=0,label=left:$\ldots$, label=right:$p+1$] (h) at (5,0) {$\otimes$};
  \path (c) edge node[below]{} (b)
          (c) edge node[below]{} (k)
        (b) edge node[below]{} (p)
           (b) edge node[below]{} (o)
         (k) edge node[below]{} (n)
        (k) edge node[below]{} (m)
        (a) edge node[below]{} (c)
        (a) edge node[below]{} (f)
      (h) edge node[below]{} (i)
      (h) edge node[below]{}  (j)
       (g) edge node[below]{} (d)
      (g) edge node[below]{} (e)
        (f) edge node[left]{} (g)
        (f) edge node[left]{} (h);
\end{tikzpicture}
       \caption{The cotree $T_{G}$}
       \label{Figu6}
\end{center}
\end{figure}
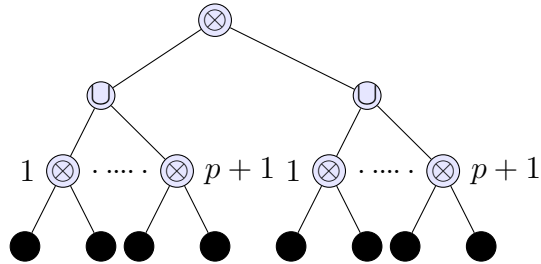

Consider now the cograph $G=(p+1)K_{2}\otimes(p+1)K_{2}$, of order $n=4p+4$, whose regular balanced cotree $T_{G}(2,p+1,0|0,0,2)$  is represented in Figure \ref{Figu6}.

\begin{Lem}
\label{lem11}
Let $G=(p+1)K_{2}\otimes(p+1)K_{2}$ be a cograph of order $n=4p+4$, for a fixed value $p\geq 1$. Then the spectrum of $G$ is
$$-(2p+1); -1; 1; 2p+3$$ with multiplicity $$1; 2(p+1); 2p; 1,$$ respectively.
\end{Lem}
{\bf Proof:} Notice that, using Theorem \ref{Tre1}, we can consider that $m=2(p+1)$, $t_{1}=\cdots=t_{m}=2$ . Then the multiplicity of $-1$ is
$$\sum_{i=1}^{m}(t_{i}-1)=2(p+1).$$
Noticing that $T_{G}(2,p+1,0|0,0,2)$ is a regular balanced cotree, we can apply Lemma \ref{lem6} to obtain that $m(1)=2p$. To obtain that $m(-(2p+1))=1$ we just execute the algorithm diagonalize with input $(T_{G},2p+1)$ and observe that it creates a single zero on the $T_G$. The eigenvalue $2p+3$ is determined by the fact that the eigenvalues must sum zero.

\begin{Thr}
Let $G=(p+1)K_{2}\otimes(p+1)K_{2}$ be the cograph of order $n=4p+4$, represented in Figure \ref{Figu6}, for a fixed value $p\geq 1$. Then $G$ is integral, non-cospectral and borderenergetic to $K_{n}$.
\end{Thr}
{\bf Proof:}
Using Lemma \ref{lem11} we have that $E(G)=8p+6$ and $E(K_{n})=2(n-1)=8p+6$. And the result follows.


\begin{thebibliography}{99}



\bibitem{Bell}
F.K. Bell, P. Rawlinson, On the multiplicities of graph eigenvalues,  {\it  Bull. Lond. Math. Soc.\/} {\bf 35} (2003) 401--408.


\bibitem{BSS2011}
T. B{\i}y{\i}ko\u{g}lu, S. K. Simi\'{c}, Z. Stani\'c, Some notes
on spectra of cographs, {\it Ars Comb.\/} {\bf 100} (2011) 421--434.


\bibitem{Bu}
C. Bu, X. Zhang, J. Zhou,
A note on the multiplicities of graph eigenvalues,
 {\it Lin. Algebra Appl. \/} {\bf 442} (2014) 69--74.



\bibitem{Chang08} Gerard. J. Chang, Liang-Hao Huang, Hong-Gwa Yeh,
On the rank of a cograph,  {\it Lin. Algebra Appl. \/} {\bf 429} (2008), 601-605.


\bibitem{Stewart}
D.G. Corneil, H. Lerchs and L. Stewart Bhirmingham;  Complement reducible
graphs, {\it Discr. Appl.  Math. \/} {\bf 3} (1981) 163--174.



\bibitem{Li2}
B. Deng, X. Li and I. Gutman,  More on  borderenergetic graphs, {\it Lin. Algebra Appl. \/}
{\bf 497} (2016) 199-208.






\bibitem{Godsil}
C. Godsil, G. Royle, {\it Algebraic Graph Theory\/}, Springer-Verlag, New York, 2001.


\bibitem{Gutman2015}
S. C. Gong, X. Li, G. H. Xu, I. Gutman, B. Furtula,
Borderenergetic graphs,
{\it MATCH Commun. Math. Comput. Chem. \/}{\bf 74} (2015) 321-332.




\bibitem{Hou}
 Y. Hou, Q. Tao, Borderenergetic threshold graphs, {\it MATCH Commun.
    Math. Comput. Chem.\/} {\bf 75} (2016) 253--262.





\bibitem{JTT2016}
D. P. Jacobs, V. Trevisan, F. C. Tura, Eigenvalue location   in
cographs, {\it  Discr. Appl. Math.\/}   (2017), $http://dx.doi.org/10.1016/ j.dam.2017.02.007.$


\bibitem{JTT2015}
D. P. Jacobs, V. Trevisan, F. Tura, Eigenvalues and energy in
    threshold graphs, {\it Lin. Algebra Appl.\/} {\bf 465} (2015)
    412--425.



\bibitem{JTT2013}
D. P. Jacobs, V. Trevisan, F. Tura, Eigenvalue location in threshold graphs,
{\it Lin. Algebra  Appl.\/} {\bf 439} (2013) 2762--2773.



\bibitem{Li}
X. Li, M. Wei, S. Gong, A computer search for the borderenergetic graphs of order 10,
{\it MATCH Commun. Math. Comput. Chem. \/}{\bf 74} (2015) 333-342.



\bibitem{Gutman2012}
X. Li, Y. Shi, I. Gutman, {\it Graph Energy\/}, Springer,
    New York, 2012.




\bibitem{Moha}
 A. Mohammadian, V. Trevisan, Some spectral properties of cographs, {\it Discr. Math.\/} {\bf 339} (2016) 1261--1264.





\bibitem{Row} P. Rowlinson,
On multiple eigenvalues of tree, {\it Lin. Algebra  Appl. \/}{\bf 432} (2010), 3007-3011.




\bibitem{Royle}
G. F. Royle;  The rank of cographs, {\it  Elec. J. Comb.\/}{ \bf 10} (2003) $N_{0}$
1. Research  paper $N$ 11:7 p.



\bibitem{Sal2015} J. Salez, Every totally real algebraic integer is a tree eigenvalue, {\it J. Comb. Theory Ser. B\/}{ \bf 111} (2015) 249–-256.

\bibitem{SF2011}
I. Sciriha, S. Farrugia,
On the spectrum of threshold graphs,
{\it ISRN Discr. Math.\/} (2011)
doi:10.5402/2011/108509.

\bibitem{Li3}
Z. Shao, F. Deng,  Correcting the number of borderenergetic graphs of order 10,
{\it MATCH Commun. Math. Comput. Chem. \/} {\bf 75} (2016) 263-266.





\end{thebibliography}
\end{document}